\documentclass{article}

\usepackage{amsmath,amsthm,amssymb,graphicx}

\newcommand{\Md}{\mathrm{maxdeg}\ }
\newcommand{\md}{\mathrm{mindeg}\ }

\newtheorem{theorem}{Theorem}
\newtheorem{definition}{Definition}
\newtheorem{conjecture}{Conjecture}

\title{A slope conjecture for links}
\author{Roland van der Veen}

\begin{document}
\maketitle
\abstract{The slope conjecture \cite{Ga11} gives a precise relation between the degree of the colored Jones
polynomial of a knot and the boundary slopes of essential surfaces in the knot complement. In this note
we propose a generalization of the slope conjecture to links. We prove the conjecture for all alternating or more generally adequate links.
We also verify the conjecture for torus links.}

\section{Introduction}

As predicted by Witten's path integral formulation \cite{Wi89}, the (colored) Jones polynomial provides many
relations between classical and quantum topology. For example the AJ conjecture relates the $q$-difference
equation for the colored Jones polynomials to the character variety of the knot group \cite{Ga04,FGL02}.
A closely related but perhaps more tractable problem is the slope conjecture \cite{Ga11}. 
The slope conjecture gives an interpretation of the growth of the maximal degree of the the colored Jones polynomials
in terms of essential surfaces in the knot complement. The purpose of this note is to suggest a generalization of the slope conjecture
to links. 

Let $L$ be a link in $S^3$ with $\ell$ components. For a link $L$ in $S^3$, let $T$ denote a tubular neighborhood of $L$ and let $M = \overline{S^3 - T}$
denote the exterior of $L$. The boundary $\partial M = \bigcup_j \partial M_j$ is a union of tori, one for each component of $L$.
Let $(\lambda_j,\mu_j)$ be the canonical meridian and longitude for the $j$-th torus component $\partial M_j$.

\begin{definition}
Suppose there exists a properly embedded essential surface $(\Sigma, \partial \Sigma) \subset (M , \partial M)$, such that every circle of
$\partial \Sigma\cap \partial M_j$ is homologous to $p_j\mu + q_j\lambda$. Then we call the sequence $p_j/q_j \in Q \cup \{\frac{1}{0}\}$
a boundary slope of $L$.
\end{definition}

Next consider the colored Jones function $J(L):\mathbb{N}^\ell \to \mathbb{Z}[v,v^{-1}]$, see Section \ref{sec.Jones} for the definition.
It sends a sequence of colors $N=(N_1,\hdots,N_\ell)$ to the unnormalized colored Jones polynomial $J_N(L;v)$. Here each component is assumed to be
$0$-framed. For example if $H$ is the Hopf link, then $\ell=2$ and 
\[J_N(H;v) = (-1)^{N_1+N_2}\frac{v^{2N_1 N_2}-v^{-2N_1 N_2}}{v^2-v^{-2}}\]
The quadratic behaviour of the maximal degree in $v$ of $J(L)$ is described in the following definition. 

\begin{definition}
\label{def.slopematrix}
A slope matrix for a function $d:\mathbb{N}^\ell \to\mathbb{N}$ is an $\ell\times \ell$ symmetric matrix $S$ if there exist infinite subsets
$U_1 \ldots U_\ell\subset \mathbb{N}$ and a constant $C$ such that for $N \in \prod_{j=1}^\ell U_j$ we have $|d(N)- N S N^t| < C|N|$. 
\end{definition}

We are interested in slope matrices for functions $\Md J(L)$. In the Hopf link example the funciton $\Md J(H)$ has a single slope matrix $S^H$ given by
$S^H = \left( \begin{array}{cc}
0 & 1 \\
1 & 0 \end{array} \right)$. Indeed, the maximal degree $J_N(H;v)$ is roughly $2N_1N_2 = N S_H N^t$ for all $N$ so we can take $U_1=U_2=\mathbb{N}$ and $C=4$. 

We are now ready to formulate the slope conjecture for links.

\begin{conjecture}{\rm (Slope conjecture for links)}\\
\label{conj.slopeLink}
If $S$ is a slope matrix for $\Md J(L)$ then there exists an essential
surface $\Sigma$ in the complement of $L$ whose boundary slope is given by the column sums
of the slope matrix. So $j$-th element of the boundary slope equals $\sum_i S_{ij}$. 
\end{conjecture}

\begin{figure}[htp]
\label{fig.Hopf}
\begin{center}
\includegraphics[width=5cm,height=5cm]{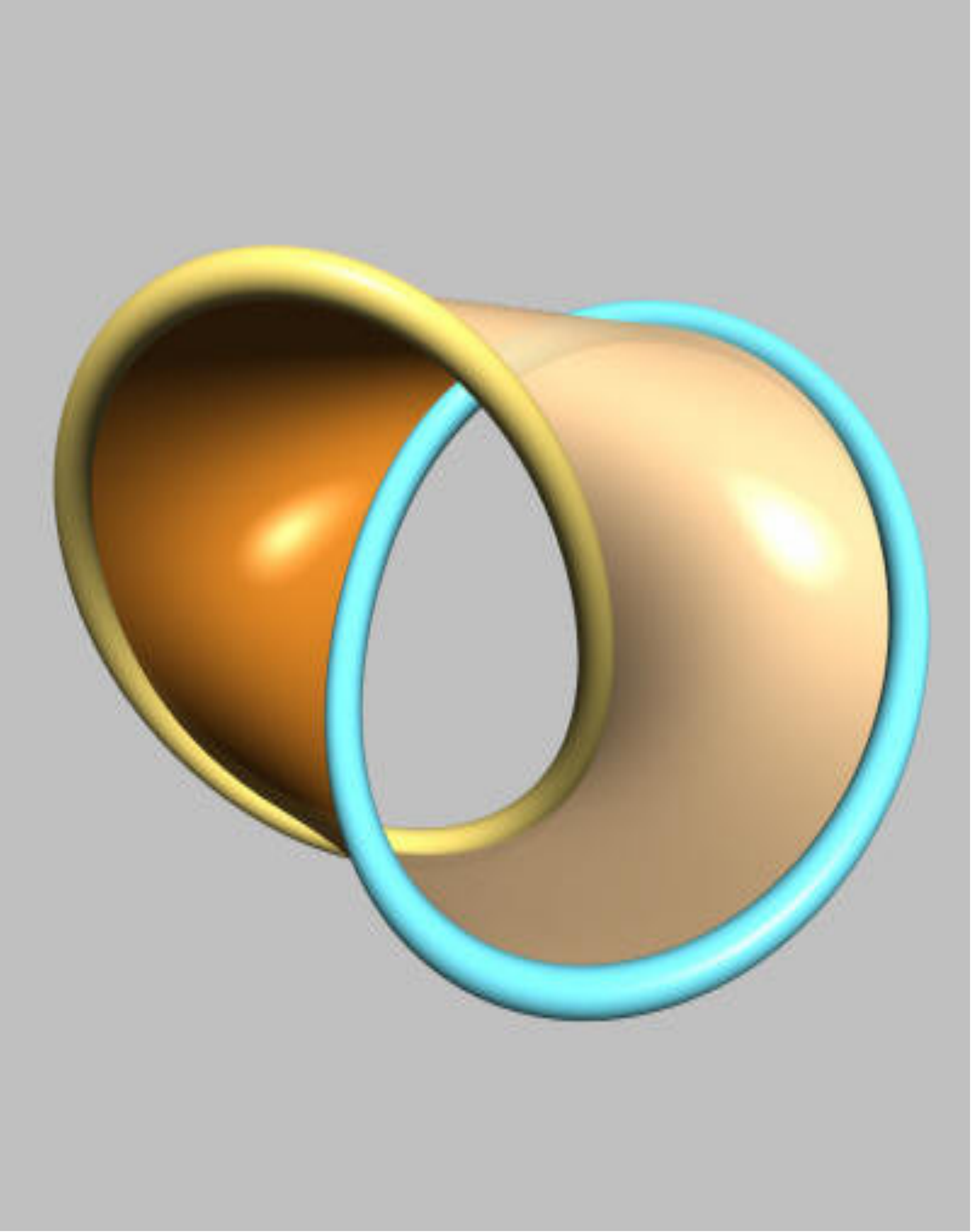}
\caption{The Hopf link together with the surface $\Sigma$ of slope $(1,1)$. Notice that a closed path on $\Sigma$ close to the green ring winds around it once.
Image created by SeifertView \cite{CW06}.}
\end{center}
\end{figure}

For example the above slope matrix $S^H$ for the maximal degree of the Hopf link should yield a surface $\Sigma$
with slope $(\frac{1}{1},\frac{1}{1})$, the column sums of the slope matrix. In Figure \ref{fig.Hopf} we 
have drawn such a surface.

In the case of a knot $\ell=1$ a slope matrix $s$ for $\Md J(K)$ 
coincides with an accummulation point of the set $\{\frac{\Md J_N(K)}{N^2}|N\in \mathbb{N}\}$. Therefore our
slope conjecture for links generalizes the original slope conjecture for knots as stated in \cite{Ga11}.
here is a similar statement for the minimal degree of the knot but this follows from the present conjecture
by replacing the knot by its mirror image. 

To provide some evidence for our conjecture we prove it for the large class of $A$-adequate links by generalizing
the proof of \cite{FKP11} to links. Minus-adequate links are a vast generalization of alternating links, see Section \ref{sec.A}
for a precise definition.

\begin{theorem}
\label{thm.Aslope}
The slope conjecture is true for all minus-adequate links.
\end{theorem}

Given the simple closed form for torus links it is not hard to check the conjecture holds in this case as well.

\begin{theorem}
\label{thm.torus}
The slope conjecture holds for all torus links.
\end{theorem}

As additional evidence for the conjecture we note that if we take two knots for which the 
slope conjecture is known then the slope conjecture also holds for their distant union. This is
because the colored Jones polynomial is multiplicative and the slopes additive. Connected sum can also be checked
along the lines of \cite{MT15}.

From these examples and theorems it may appear that the degree is rather well behaved and easy to compute from the state sums
and the surfaces are simply some Seifert-type state surfaces. In the general case this is quite far from the truth. 
Massive cancellations will usually render the ordinary state sums useless. Also the surfaces are in general quite complicated as the slopes are often
rational numbers with high denominator \cite{GV14}. As such the slope conjecture presents a major challenge to both classical and quantum topology.
\\

\noindent
{\bf Acknowledgement} A preliminary version of this work was presented at the Moab topology conference 2012. 
I would like to thank the organizers Nathan Geer and Jessica Purcell for their hospitality and for providing an
inspiring atmosphere. I would also like to thank Oliver Dasbach, Dave Futer, Stavros Garoufalidis, 
Effie Kalfagianni and Jessica Purcell for stimulating conversations. This work was supported by the Netherlands Organisation for
Scientific Research.

\section{Colored Jones polynomial for links}
\label{sec.Jones}

Since multiple conventions are in use in the literature we briefly give a definition of the colored Jones polynomial of a link.
The present definition is perhaps not the most natural but has the advantage of being relatively concrete and easy to state.
Assuming the reader is familiar with simple skein theory \cite{Li97} it is easiest to construct the 
colored Jones polynomial from the colored Kauffman bracket as follows.

Recall that the skein space of the annulus can be regarded as the polynomial ring in the variable $x$
corresponding to its core. Given a framed knot diagram $D$ we can insert an element $P(x)$ of the skein of
the annulus and take the Kauffman bracket in the plane. The resulting Laurent polynomial in $A$ will be written 
$\langle D(P(x))\rangle$. More generally for a $\ell$-component link diagram $D$ 
and polynomial $P_j(x)$, on for each component we will write
$\langle D(P_1(x),\ldots,P_\ell(x))\rangle(A)$.

Because of representation theory a special role is played by the Chebyshev polynomials $S_n(x)$ defined by 
$S_{n+1}(x) = x S_n- S_{n-1}$ and $S_0 = 1$ and $S_1 = x$.

\begin{definition}
\begin{enumerate}
Let $L$ be an $\ell$-component link in $S^3$ with diagram $D$ where the writhe of the $j$-th component is $w_j$.
\item The $N = (N_1,\ldots, N_\ell)$-colored Kauffman bracket of $L$ is defined to be
\[
\langle L \rangle_N(A) = \langle D(S_{N_1},\ldots S_{N_\ell}) \rangle(A)
\]
\item Set $v=A^{-1}$. The $N$-colored Jones polynomial of an $\ell$-component link $L$ is defined to be
\[
J_N(L;v) = \prod_{j=1}^\ell((-1)^{N_j-1}v^{N_j^2-1})^{w_j}\langle L \rangle_{N_1-1,N_2-1,\ldots N_\ell-1}(v)
\]
\end{enumerate}
\end{definition}

Note the shift from $N$ to $N-1$. The version of the colored Jones polynomial we defined is unnormalized and framing independent because of the writhe term.
The value of the $N$-colored unknot is $\frac{v^{2N}-v^{-2N}}{v^2-v^{-2}}$.
It is perhaps more standard to replace $A$ by $q^{-\frac{1}{4}}$ but the variable $v=A^{-1}$ neatly absorbs
the factor of $4$ otherwise present in the slope conjecture.

\section{Proof of the slope conjecture for minus-adequate links}
\label{sec.A}
Recall that in computing the Kauffman bracket each crossing gets resolved in two different ways:
the plus-resolution and the minus-resolution (also known as the $A$- and $B$-resolutions). 
The result is a linear combination, the state sum, where each $\pm$ resolution contributes $A^{\pm 1}$. 
A diagram is called minus-adequate (also $B$-adequate) if after taking the minus-resolution everywhere, 
the two arcs replacing each crossing always belong to different circles. In other words,
the all minus-resolution results in the maximal number of state circles. As such it must contribute 
the lowest power of $A$, namely $A^{-c-2p}$ where $c$ is the number of crossings and
$p$ is the number of circles in the state. See Lemma 5.4 \cite{Li97}.

Finally recall that the minus-state surface is the surface obtained from the diagram by the following steps.
First do the all minus-resolution and turn each of the resulting circles into disks. Next at each former crossing,
attach a half twisted band connecting a pair of disks. This surface is essential 
according to Ozawa \cite{Oz11}, see also \cite{FKP13} Theorem 3.19.

We are now ready to prove the slope conjecture for minus-adequate links. Our proof generalizes that for knots
found in \cite{FKP11}. 

Let $D$ be a $-$-adequate link diagram whose all minus-resolution gives
rise to $p$ state circles. We assume the link has $\ell$ components and we denote the diagram of the $j$-th component by $D_j$.
Also $c^+_{ij}$ and $c^-_{ij}$ are the number of positive and negative crossings between $D_i$ and $D_j$. 
The total number of such crossings is $c_{ij}$ and the writhe is $w_j = c^+_{jj}-c^-_{jj}$.

If we denote by $D^N$ the $N$-parallel of the diagram (component $j$ gets replaced by $N_j$ parallel copies)
then 

\[\Md J_{N}(L;v) = \sum_{j}w_j(N_j^2-1)-\md\langle D^{N-1} \rangle(A)\]

This is because we can expand the Chebyshev polynomials into a linear combination of minus-adequate diagrams.
Since the all minus-resolution yields the lowest degree term for each diagram, the overall lowest degree is that
of the leading term of the Chebyshev polynomial. Moreover we know exactly what this minimal degree of the Kauffman bracket is:
minus the number of crossings minus twice the number of state circles:

\[-\md\langle D^{N-1} \rangle(A) = c+2p = \sum_{i\leq j}(N_i-1)(N_j-1)c_{ij}-2p\]

We can estimate $0\leq 2p \leq 2p'\mathrm{max}_jN_j$ where $p'$ is the number of state circles in the minus-resolution of $D$ itself.
 
Summing up we found the following estimate for the maximal degree in $v$ of $J_N(L)$:

\[
\Md J_N(L) = 2\sum_{i}N_i^2c^+_{ii}+\sum_{i< j}N_iN_jc_{ij} + \mathcal{O}(|N|)
\]

In other words the matrix $S$ defined below is a slope matrix in the sense of Definition \ref{def.slopematrix}. 
\[S_{ij} =  \begin{cases}
      \frac{1}{2}c_{ij} & \text{if $i \neq j$}\\
      2c^+_{ii} & \text{if $i = j$}
   \end{cases}
\]

For this we take the infinite subsets $U_j$ from the definition to be $\mathbb{N}$ and choose some constant $C>2p'$.

In order to prove Theorem \ref{thm.Aslope} we now need to find an essential surface
whose slope at the $j$-th component equals $\sum_{i}S_{ij} = \frac{1}{2}\sum_{i\neq j}c_{ij}+ 2c^+_{jj}$.
For this we can use the all minus-state surface introduced above.
The slope of this surface is found by calculating the linking number between the $j$-th component $D_j$ and
a curve following it along the surface. For every negative crossing between $D_j$ and itself we find a contribution of
$2$. For positive such crossings we get a contribution of $0$. For any type of crossing between $D_j$ and another component
$D_i$ we obtain a contribution $\frac{1}{2}$.
In total one thus gets a slope of $2c^-_{jj}+ \frac{1}{2}\sum_{i\neq j}c_{ij}$ as required.

\section{The case for torus links}

We consider the $(r,s)$-torus link $T^r_s$ with $r,s \in \mathbb{Z}$ and $s\geq1$ to be the closure of the braid
$(\sigma_1\sigma_2\cdots \sigma_{s-1})^r$ in $S^3$. In case $r<0$ the standard braid diagram is minus-adequate.
If follows from the formula of the Jones polynomial below that the case $r>0$ is not minus-adequate and hence
provides additional evidence to our version of the slope conjecture for links.

The link $T^r_s$ has $g = \gcd(r,s)>0$ components. If we define $a = r/g$ and $b = s/g$ 
then the colored Jones polynomial colored by $N = (N_1,N_2,\cdots,N_g)$ is given by the formula \cite{vV08}

\[J_N(T^r_s) = v^{ab(|N|^2-g)}\sum^{(|N|-g)}_{k = -(|N|-g)}{g\choose k}_N v^{-ak(bk+2)}[bk+1]\]
where the summation variable $k$ takes steps of two and we used the notation $|N| = N_1+N_2+\hdots N_g$ and $|N|^2 = N_1^2+N_2^2+\hdots N_g^2$
and finally ${g \choose k}_N$ is the coefficient of $v^{2k}$ in $\prod_{j = 1}^g [N_j]$.

For $r<0$ we see that the maximal degree is given by the term $k = (|N|-g)$. This term has maxdegree
$-ab(|N|^2-g)+ab(|N|-g)^2+ 2a(|N|-g)+ 2b(|N|-g) = ab\sum_{i,j}N_iN_j +\mathcal{O}(|N|)$. Hence the unique slope matrix $S$ is given
by $S_{ij} = ab (1-\delta_{ij})$. This is in agreement with Theorem \ref{thm.Aslope} applied to the standard braid diagram.
Hence the state-surface provides the surface with the correct slopes $ab(g-1)$.

The case $r>0$ is more interesting. The formula for the maximal degree is not a quadratic polynomial in $N$ but
is piecewise quadratic depending on the parity of $|N|-g$. This already shows that such links cannot be minus-adequate 
since then the degree would be a quadratic polynomial in $N$.

More specifically if $|N|-g$ is even then maximal degree is given by the term $k = 0$ and equals $ab(|N|^2-g)$.
Care has to be taken with the case $b=1$ since there the term $k=-2$ contributes with the same degree but luckily 
with a different coefficient: For $g>1$ and $N_j>1$ we have ${g\choose 0}_N\neq {g\choose -2}_N$.

Next we look at the case $|N|-g$ odd. The term $k = -1$ contributes the maximal degree which now equals
$ab(|N|^2-g) + ab-2a+2b-4$, except in the cases $ab\leq 2$. By the symmetry between $a$ and $b$ we may assume that
$b = 1$ in which case the $k = -1$ term vanishes. The $k = 1$ term then contributes $ab(|N|^2-g) - ab-2a+2b$

In either case the correct slope matrix is $S = ab I$. To complete the slope conjecture we need to find 
an incompressible surface bounding the link whose boundary slopes are all $ab$. Consider the canonical
annuli defined by taking the complement of the $(r,s)$ torus link viewed as sitting on the torus surface. 
This surface is certainly essential and each component is seen to have slope $ab$ as required. 
This proves Theorem \ref{thm.torus}.

With a bit more work one may be able to prove the slope conjecture for all zero-volume links as
these are very similar to torus links. Their colored Jones polynomials can be readily computed by repeatedly
cabling and taking connected sums, see \cite{vV08}.

\section{Further directions}

Computing boundary slopes of surfaces in link complements is an interesting problem in classical topology. 
Outside the knot case not much is known except for some work on the two-bridge link case \cite{HS07}. 
Through the slope conjecture one can use the colored Jones polynomial to explore such slopes further.
At least in some cases \cite{GV14} the Jones polynomial suggests the existence of some very non-trivial slopes,
thus providing new challenges to classical topology.

One wonders if there a way to extract these slopes from the tropical geometry of the A-ideal that replaces 
the A-polynomial in the link case. This is natural since one expects an AJ-type conjecture for links.

On the colored Jones side an important question is whether the degree is still a (multivariate) quadratic quasipolynomial
as in the knot case \cite{Ga113}. At least for minus-adequate links and the torus links we considered this is true. 
Such a question brings us back to a possible link version of the AJ conjecture. An important first step towards such a conjecture
has been set for two-bridge links in \cite{LT11}.

It would also be interesting to go beyond the slope conjecture as in \cite{GL11} and consider stabilization
properties for links. At least for alternating links their (multivariate) heads and tails and beyond 
can readily be computed.

Finally the question of the behaviour of the degree can be posed for any quantum invariant, not just the colored Jones polynomial.
It would be interesting to see what happens for other knot polynomials and their categorifications. A first approach to this question
in the colored HOMFLY case has been made in \cite{vV15}.

\bibliography{biblio}{}
\bibliographystyle{plain}

\end{document}